\documentclass[12pt]{article}
\usepackage{enumerate}
\usepackage{amsmath,amssymb,amsthm, mathrsfs, mathtools}
\usepackage{latexsym}
\usepackage{graphics,graphicx}
\usepackage{hyperref}
\usepackage[bf, small]{titlesec}
\usepackage{authblk} 
\usepackage{soul}
\usepackage{kotex}
\usepackage{lineno}

\theoremstyle{plain}
\newtheorem{Thm}{Theorem}[section]

\newtheorem{Lem}[Thm]{Lemma}
\newtheorem{Prop}[Thm]{Proposition}
\newtheorem{Cor}[Thm]{Corollary}

\newtheorem{Conj}[Thm]{Conjecture}

\theoremstyle{definition}

\usepackage{standalone}
\usepackage{tikz}
\usetikzlibrary{arrows,positioning,matrix,fit,backgrounds,shapes,shapes.geometric, intersections}

\usetikzlibrary{shapes.callouts,decorations.pathmorphing} 
\usepackage{calc} 
\usetikzlibrary{decorations.markings} 
\tikzstyle{vertex}=[circle, draw, inner sep=0pt, minimum size=6pt] 

\usetikzlibrary{arrows,matrix} 
\newtheorem*{claim1}{Claim A}

\newcommand{\olr}[1]{\overleftrightarrow{#1}}

\title{Digraphs in which every $t$ vertices have exactly $\lambda$ common out-neighbors}

\author[1]{Myungho Choi}
\author[2]{Hojin Chu}
\author[1]{Suh-Ryung Kim}

\affil[1]{\footnotesize Department of Mathematics Education, Seoul National University, Seoul 08826, Rep. of Korea}
\affil[ ]{\footnotesize\textit{nums8080@snu.ac.kr, srkim@snu.ac.kr}}
\date{}

\affil[2]{\footnotesize School of Computational Sciences, Korea Institute for Advanced Study (KIAS), Seoul, Rep. of Korea}
\affil[ ]{\footnotesize\textit{hojinchu@kias.re.kr}}

\definecolor{LemonChiffon}{rgb}{100, 98, 80}
\definecolor{myblue}{rgb}{0,0.4,0.8}
\definecolor{orange}{rgb}{1, 0.4, 0}
\definecolor{mygreen}{rgb}{0, 0.8, 0}
\definecolor{myred}{rgb}{204, 0, 0}
\definecolor{violet}{RGB}{0.4,0.2,1}
\definecolor{brown}{rgb}{0.6, 0.4, 0}

\newcounter{statement}

\begin{document}
\maketitle
\begin{abstract}
We say that a digraph is a $(t,\lambda)$-liking digraph if every $t$ vertices have exactly $\lambda$ common out-neighbors.
In 1975, Plesn\'{i}k [Graphs with a homogeneity, 1975. {\it Glasnik Mathematicki} 10:9-23] 
proved that any $(t,1)$-liking digraph is the complete digraph on $t+1$ vertices for each $t\geq 3$.
Choi {\it et al}. [A digraph version of the Friendship Theorem, 2025. {\it Discrete mathematics}, 348(1), 114238]
showed that a $(2,1)$-liking digraph is a fancy wheel digraph or a $k$-diregular digraph for some positive integer $k$.
In this paper, we extend these results by completely characterizing the $(t,\lambda)$-liking digraphs with $t \geq \lambda+2$ and giving some equivalent conditions for a $(t,\lambda)$-liking digraph being a complete digraph on $t+\lambda$ vertices.
 \end{abstract}



\section{Introduction}
In this paper, for graph-theoretical terminology and notations not defined, we follow \cite{bang2018classes} and \cite{bondy2010graph}.
Neither graphs nor digraphs in this paper have loops, multiple edges, or multiple arcs.

This paper was motivated by the Friendship Theorem introduced by Erd\"{o}s {\it et al.} \cite{erd1966r} in 1966.
Using graph-theoretical terminology, the Friendship Theorem 
can be described as follows: if any pair of vertices in a graph has exactly one common neighbor, then there exists a vertex adjacent to the others.
A graph that satisfies the hypothetical part of the Friendship Theorem is called a ``friendship graph", that is,
a {\it friendship graph} is a graph such that every pair of vertices has exactly one common neighbor. 
Many  variants  of a friendship graph  have  been  studied (see \cite{bose1969graphs,CARSTENS1977286,DELORME1984261,FERNAU2008734,GUNDERSON2016125,hartke2008question,LI20121892,muller1974strongly,plesnik1975graphs,sudolsky1978generalization}). 
 
A digraph version of the Friendship Theorem was introduced and studied as a variant (\cite{choi2023digraph, muller1974strongly, plesnik1975graphs}).
A digraph is called a
{\it $(t,\lambda)$-liking digraph} if every $t$ vertices have exactly $\lambda$ common out-neighbors for some positive integers $t$ and $\lambda$.
In 1974,
 M\"{u}ller and Pelant~\cite{muller1974strongly} studied the case where tournaments are $(t,\lambda)$-liking digraphs and showed the non-existence of such tournament for $t\geq 3$.
In 1975, Plesn\'{i}k \cite{plesnik1975graphs} characterized the $(t,1)$-liking digraphs for each integer $t\ge 3$ by proving that every $(t,1)$-liking digraph with $t\ge 3$ is the complete digraph on $t+1$ vertices, which is a generalization of the result given by M\"{u}ller and Pelant.
He proposed as an open question to take care of the case $t=2$.
In 2025, Choi {\it et al.} \cite{choi2023digraph} settled  the question by completely characterizing the $(2,1)$-liking digraphs.
   
In this paper, we extend the results of Plesn\'{i}k and Choi {\it et al.} by completely characterizing the $(t,\lambda)$-liking digraphs for any positive integers $t$ and $\lambda$ satisfying $t\ge \lambda+2$ as follows.
   \begin{Thm}\label{thm:main1}
If $t\ge \lambda+2$, then the complete digraph on $t+\lambda$ vertices is the only $(t,\lambda)$-liking digraph.
\end{Thm}
Indeed, we give some equivalent conditions for a $(t,\lambda)$-liking digraph being the complete digraph on $t+\lambda$ vertices.

\begin{Thm}\label{thm:main2}
Let $D$ be a $(t,\lambda)$-liking digraph for some positive integers $t, \lambda$ with $t\ge 2$.
Then the following are equivalent.
\begin{enumerate}[(a)]
	\item $D$ is complete on $t+\lambda$ vertices, that is, $D \cong \olr{K}_{t+\lambda}$.
	\item $D$ is a $(t-1,\lambda+1)$-liking digraph.
	\item $d^+(v)=t+\lambda-1$ for each vertex $v$.
\end{enumerate}
Furthermore, (a) is equivalent to the condition (d) if $(t, \lambda) \neq (2,1)$; (e) if $t\ge 3$, where
\begin{enumerate}[(a)]
	\item[(d)] there is a vertex $v$ satisfying $N^+(v)=V(D)\setminus \{v\}$;
	\item[(e)] $D$ is diregular.
\end{enumerate}

\end{Thm} 
We will prove the above two theorems in Section~\ref{sec:proofs}.

We note that the equivalence of (a) and (e) in Theorem~\ref{thm:main2} implies that if $t\ge 3$, then the complete digraph on $t+\lambda$ vertices is the only diregular $(t,\lambda)$-liking digraph.
This reveals a key structural rigidity specific to the directed setting: diregularity alone forces completeness when $t\ge 3$, a phenomenon that breaks down at $t=2$ (see Figure~\ref{fig:3-regular}).

There is a variant of a friendship graph, called a ``generalized friendship graph".
A {\it generalized friendship graph} is a graph if any $t$ vertices have exactly $\lambda$ common neighbors for some positive integers $t$ and $\lambda$. 
It can be said that a $(t, \lambda)$-liking digraph is a digraph version of a generalized friendship graph.
Considering this observation,
we refer this graph as a $(t,\lambda)$-friendship graph to make the concept clearer.
The structure of $(t,\lambda)$-friendship graph is well-known.
In \cite{CARSTENS1977286} and \cite{sudolsky1978generalization}, it was shown that a $(t, \lambda)$-friendship graph is the complete graph on $t+\lambda$ vertices if $t\ge 3$. 
Theorem~\ref{thm:main1} shows that the digraph version of this theorem is true if $t\ge \lambda+2$.
Moreover, the equivalence of (a) and (e) in Theorem~\ref{thm:main2} generalizes the theorem for $(t,\lambda)$-friendship graphs since a $(t,\lambda)$-friendship graph becomes a diregular $(t,\lambda)$-liking digraph by replacing each edge with a directed cycle of length two.

	\section{Preliminaries}\label{sec:pre}

This section describes the existing results needed to prove Theorems~\ref{thm:main1} and \ref{thm:main2} and the results that can be simply derived from them.

\begin{Lem}[Lemma~6.2 in \cite{plesnik1975graphs}]\label{lem:outdegree}
	Let $D$ be a $(t,\lambda)$-liking digraph.
	Then $|V(D)|\ge t+\lambda$ and $\delta^+(D) \ge t+\lambda-1$.
\end{Lem}

\begin{Lem}[Lemmas~6.4 and 6.6 in \cite{plesnik1975graphs}]\label{lem:inequalities}
		Let $D$ be a $(t,\lambda)$-liking digraph.
	Then the following hold. 
	\begin{enumerate}[(1)]
		\item \[\sum_{v\in V(D)} {d^-(v) \choose t} = \lambda {|V(D)| \choose t}. \]
		\item For each $v\in V(D)$, 
		\[{d^-(v) \choose t-1} \le {d^+(v) \choose \lambda}. \]
	\end{enumerate}
\end{Lem}

\begin{Prop}[Lemma~6.3 in \cite{plesnik1975graphs}]\label{prop:jan}
	Let $D$ be a $(t,\lambda)$-liking digraph.
	Then for any $1\le i < t$ and $\{v_1, \ldots, v_i\} \subseteq V(D)$, the subdigraph of $D$ induced by the common out-neighbors of $v_1, \ldots, v_i$ is a $(t-i,\lambda)$-liking digraph.
\end{Prop}

We make a meaningful observation from Lemma~\ref{lem:outdegree} and Proposition~\ref{prop:jan}.

\begin{Prop}\label{prop:out-degree}
Let $D$ be a $(t,\lambda)$-liking digraph.
For each $1\leq i < t$, every $t-i$ vertices have at least $\lambda+i$ common out-neighbors in $D$.
\end{Prop}

\begin{proof}
	Fix $i \in \{1,\ldots, t-1\}$.
	Take $t-i$ vertices.
	Then the subdigraph of $D$ induced by the common out-neighbors of the $t-i$ vertices is a $(i,\lambda)$-liking digraph by Proposition~\ref{prop:jan}, and so it has at least $\lambda+i$ vertices by Lemma~\ref{lem:outdegree}. 
	Thus every $t-i$ vertices have at least $\lambda+i$ common out-neighbors in $D$.
\end{proof}

Given a digraph $D$ and vertices $u$ and $v$ of $D$,
we use the expression $u \to v$ when $(u,v) \in A(D)$. 
When representing negation, add a slash ($/$) to the symbol.

 A {\it $k$-diregular digraph} is a digraph in which each vertex has outdegree $k$ and indegree $k$ for a positive integer $k$.
	A digraph is said to be {\it diregular} if it is a $k$-diregular digraph for some positive integer $k$.

Let $v$, $k$, $\lambda$, and $t$ be positive integers such that $v>k\ge t$.
A {\it $t$-$(v,k,\lambda)$ design} $\mathfrak{D}=(X,\mathcal{B})$ consists of a set $X$ of elements, called {\it varieties}, and a collection $\mathcal{B}$ of subsets of $X$, called {\it blocks}, such that the following conditions are satisfied:
	\begin{itemize}
	\item $|X|=v$;
	\item Each block consists of exactly the same number $k$ of varieties;
	\item Every set of $t$ distinct varieties appears simultaneously in exactly the same number $\lambda$ of blocks.
	\end{itemize}
The general term {\it $t$-design} is used to indicate any $t$-$(v,k,\lambda)$ design.

It is a well-known fact that the number of blocks in a $t$-$(v,k,\lambda)$ design is
$\lambda {v \choose t}/{k\choose t}$ and each variety appears exactly the same number $\lambda {v-1 \choose t-1}/{k-1 \choose t-1}$ of blocks.
A $t$-$(v,k,\lambda)$ design is said to be {\it symmetric} if the number of its blocks is equal to the number of its varieties.
That is, \[v=\frac{\lambda {v \choose t}}{{k\choose t}}. \]
Then each variety of a symmetric $t$-$(v,k,\lambda)$ appears in exactly $k$ blocks.  \cite{stinson2004combinatorial}

\begin{Lem}\label{lem:design}
	If a $k$-diregular $(t, \lambda)$-liking digraph of order $n$ exists, then there is a symmetric $t$-$(n,k,\lambda)$ design.
\end{Lem}
\begin{proof}
	Suppose that there is a $k$-diregular $(t, \lambda)$-liking digraph $D$ of order $n$.
	To form a symmetric $t$-$(n,k,\lambda)$ design, we take the vertices of $D$ as varieties and the in-neighborhoods of the vertices as blocks.
	Since $D$ is $k$-diregular, each block has size $k$.
	Moreover, since $D$ is a $(t, \lambda)$-liking digraph, each $t$-subset of varieties appears in exactly $\lambda$ blocks.
	Thus we have obtained a symmetric $t$-$(n,k,\lambda)$ design.
\end{proof}

\begin{Prop}[Lemma~2.6 in \cite{hughes1965t}]\label{prop:trivial}
	If a $t$-$(v,k,\lambda)$ design with $t\ge 3$ is symmetric, then $k\ge v-1$.
	\end{Prop}

A digraph $D$ is {\it complete} if, for every pair $x$, $y$ of distinct vertices of $D$, both $x \to y$ and $y \to x$.
The complete digraph on $n$ vertices is denoted by $\overleftrightarrow{K}_n$.

By Lemma~\ref{lem:design} and Proposition~\ref{prop:trivial}, the following holds.
\begin{Prop}\label{prop:diregular}
	Let $D$ be a diregular $(t,\lambda)$-liking digraph with $t\ge 3$.
	Then $D$ is the complete digraph on $t+\lambda$ vertices.
\end{Prop}
\begin{proof}
	Suppose that $D$ is $k$-diregular for some positive integer $k$.
	Then, since $D$ has neither loops nor multiple arcs, $k < |V(D)|$.
	By Lemma~\ref{lem:design} and Proposition~\ref{prop:trivial}, $k \ge |V(D)|-1$.
	Thus $k=|V(D)|-1$.
	That is, $N^+(v)=V(D)-\{v\}$ for each vertex $v$.
	Therefore $D$ is complete.
	Since $D$ is a $(t,\lambda)$-liking digraph, $|V(D)|=t+\lambda$.
\end{proof}

\section{Proofs of Theorems~\ref{thm:main1} and \ref{thm:main2}}\label{sec:proofs}

The case where $t=1$ for both Theorems~\ref{thm:main1} and \ref{thm:main2} is excluded, so we only consider $(t, \lambda)$-liking digraphs with $t\ge 2$ and $\lambda\ge 1$.

\begin{Lem}\label{lem:tool}
	Let $D$ be a $(t,\lambda)$-liking digraph.
	If $d^+(v)\le d^-(v)$ for a vertex $v$, then \[0\le (d^+(v)-t-\lambda+1)(\lambda-t+1).\]
	Furthermore, if the conditional inequality is strict, then the resulting inequality is also strict.
	\end{Lem}
\begin{proof}
	Suppose that a vertex $v$ satisfies	
	\[d^+(v)\le d^-(v).\]
	Then
	\[ {d^+(v) \choose t-1}\le{d^-(v) \choose t-1}.\]
	Thus, by Lemma~\ref{lem:inequalities}(2),
	\[ {d^+(v) \choose t-1} \le {d^+(v) \choose \lambda} \]
	and so 
	\[\left|\frac{d^+(v)}{2}-\lambda\right| \le \left|\frac{d^+(v)}{2}-t+1\right|.\]
	Therefore 
	\[0 \le \left(\frac{d^+(v)}{2}-t+1\right)^2-\left(\frac{d^+(v)}{2}-\lambda\right)^2=(d^+(v)-t-\lambda+1)(\lambda-t+1).  \]
	The ``furthermore" part is obvious.
\end{proof}

\begin{Lem}\label{lem:eulerian}
	Let $D$ be a $(t,\lambda)$-liking digraph.
	If $t\ge \lambda+1$ or $d^+(v)=t+\lambda-1$ for each vertex $v$, then $d^+(v)=d^-(v)$ for every vertex $v$ in $D$.
\end{Lem}
\begin{proof}
	Assume $t \ge \lambda+1$.
	To the contrary, suppose that there is a vertex whose outdegree and indegree are different.
	Then there is a vertex $v$ such that $d^+(v)<d^-(v)$.
Thus 
	\[0< (d^+(v)-t-\lambda+1)(\lambda-t+1)  \]
	by the ``furthermore" part of Lemma~\ref{lem:tool}. 
	Then $d^+(v) \neq t+\lambda-1$.
	Since $d^+(v) \ge t+\lambda-1$ by Lemma~\ref{lem:outdegree}, $\lambda+1 >t$.
\end{proof}

\begin{Lem}\label{lem:complete,outdegree}
	Let $D$ be a $(t,\lambda)$-liking digraph.
	If $d^+(v)=t+\lambda-1$ for each vertex $v$, then $D$ is the complete digraph on $t+\lambda$ vertices.
\end{Lem}
\begin{proof}
	Suppose that $d^+(v)=t+\lambda-1$ for each vertex $v$.
	Then $d^-(v)=t+\lambda-1$ for each vertex $v$ by Lemma~\ref{lem:eulerian}.
	Therefore \[\sum_{v\in V(D)} {d^-(v) \choose t} = \sum_{v\in V(D)} {t+\lambda-1 \choose t} =|V(D)|{t+\lambda-1 \choose t}. \]
	Now, by Lemma~\ref{lem:inequalities}(1),
	\[\frac{t}{|V(D)|}{|V(D)| \choose t}=\frac{t}{\lambda}{t+\lambda-1 \choose t}.\]
	We note that 
	\[ {|V(D)|-1 \choose t-1}= \frac{t}{|V(D)|}{|V(D)| \choose t} \quad \text{and} \quad \frac{t}{\lambda}{t+\lambda-1 \choose t}={t+\lambda-1 \choose t-1}.\]
	Thus \[ {|V(D)|-1 \choose t-1}={t+\lambda-1 \choose t-1}.\]
	and so  $|V(D)|=t+\lambda$.
	Since $d^+(v)=t+\lambda-1$ for each vertex $v$, $N^+(v)=V(D)\setminus\{v\}$.
	Therefore $D \cong \olr{K}_{t+\lambda}$.
	\end{proof}

Now, we prove Theorem~\ref{thm:main1}.

\begin{proof}[Proof of Theorem~\ref{thm:main1}]
	Suppose that $D$ is a $(t,\lambda)$-liking digraph with $t\ge \lambda+2$ and take a vertex $v$.
	By Lemma~\ref{lem:eulerian},  
	\[d^+(v)=d^-(v). \]
	Then, by Lemma~\ref{lem:tool},
	\[0\le (d^+(v)-t-\lambda+1)(\lambda-t+1).\] 
	Since $t\ge \lambda+2$ by the hypothesis, $\lambda-t+1<0$ and so
	$d^+(v) \le  t+\lambda-1$.
	Thus, by Lemma~\ref{lem:outdegree}, $d^+(v) = t+\lambda-1$.
	Since $v$ was arbitrarily chosen, $d^+(u) = t+\lambda-1$ for each vertex $u$.
	Therefore $D$ is isomorphic to $\olr{K}_{t+\lambda}$ by Lemma~\ref{lem:complete,outdegree}.
\end{proof}

To prove Theorem~\ref{thm:main2}, we go further to derive more results.
For a vertex set $S$ in a digraph $D$, we denote the set of common out-neighbors (resp.\ in-neighbors) of $S$ in $D$ by $CN_D^+(S)$ (resp.\ $CN_D^-(S)$), that is, \[CN_D^+(S):= \bigcap_{u \in S} N_D^+(u) \quad \text{and} \quad CN_D^-(S):= \bigcap_{u \in S} N_D^-(u).\]
If no confusion is likely, we write $CN^+(S)$ and $CN^-(S)$ instead of $CN_D^+(S)$ and $CN_D^-(S)$, respectively.

\begin{Lem}\label{lem:complete,twoliking}
	Let $D$ be a $(t,\lambda)$-liking digraph.
	If $D$ is a $(t-1,\lambda+1)$-liking digraph, then $D$ is the complete digraph on $t+\lambda$ vertices. 
\end{Lem}
\begin{proof}
Suppose that $D$ is a $(t-1,\lambda+1)$-liking digraph.
We first claim the following.
\begin{claim1}
Let $S$ be a set of $t-1$ vertices in $D$.
Then the subdigraph of $D$ induced by the common out-neighbors of the vertices in $S$ is complete.
\end{claim1}
\begin{proof}[Proof of Claim A]
Since $D$ is a $(t-1,\lambda+1)$-liking digraph and $|S|=t-1$, 
\begin{equation}\label{eq:1}
	\left|CN^+(S)\right|=\lambda+1. 
\end{equation}
We take a common out-neighbor $v$ of the vertices in $S$.
Then $v \not\in S$ and so  $|S \cup \{v\}|=t$.
Thus, since $D$ is a $(t,\lambda)$-liking digraph,
\[\left| CN^+(S) \cap N^+(v) \right| =\lambda. \]
By \eqref{eq:1}, 
\[CN^+(S) \cap N^+(v)=CN^+(S) -\{v\}.\]
Since $CN^+(S) \cap N^+(v) \subseteq  N^+(v)$,
we have
\[CN^+(S) -\{v\} \subseteq N^+(v).\]
Therefore the vertices in $CN^+(S)$ other than $v$ are out-neighbors of $v$.
Since $v$ was arbitrarily chosen from $CN^+(S)$, the subdigraph of $D$ induced by $CN^+(S)$ is complete.
\end{proof}

Let $U$ be a maximum subset of $V(D)$ inducing a complete subdigraph in $D$.
It suffices to show that $U=V(D)$.
By Claim~A, $|U|\ge \lambda+1$.
We may assume $t \le \lambda+1$, or else we are done by Theorem~\ref{thm:main1}.
Thus $|U| \ge t$.

To the contrary, suppose that there is a subset $S$ of $U$ such that $|S|=t-1$ and $CN^+(S)\not \subseteq U$.
Take \[x \in CN^+(S)-U.\]
By the choice of $U$, $U-S \subseteq CN^+(S)$ and so $x \in CN^+(U-S)$ by Claim~A.
Since $x\in CN^+(S)$, $x \in CN^+(U)$.
Now, we take $y \in U$.
Since $|U|\ge t$, there is a subset $S'$ of $U-\{y\}$ such that $|S'|=t-1$.
Then, by the choice of $U$, $U-S' \subseteq CN^+(S')$. 
Since $x \in CN^+(U)$, $x\in CN^+(S')$ and so, by Claim~A, $x \in CN^-(U-S')$.
Thus $x \in CN^-(U-S') \subseteq N^-(y)$.
Since $y$ was arbitrarily chosen in $U$, $x \in CN^-(U)$.
Then $U \cup \{x\}$ induces a complete subdigraph in $D$, which contradicts the choice of $U$.
Therefore we may conclude that for every subset $S$ of $U$ with $|S|=t-1$, \[CN^+(S) \subseteq U.\]
Hence $|U|=t+\lambda$ since $D$ is a $(t-1,\lambda+1)$-liking digraph.

Suppose for contradiction that there is a vertex $x$ in $V(D)-U$.
Take a subset $S \in U$ with $|S|=t-1$.
Then $CN^+(S)=U-S$ and so $CN^+(S\cup\{x\}) \subseteq U-S$.
Since $|S\cup \{x\}|=t$, $|CN^+(S\cup\{x\})|=\lambda$.
Since $|U-S|=\lambda+1$, there is a unique vertex in $U-S$ that is not an out-neighbor of $x$.
If there were at least two vertices $y$ and $z$ in $U$ such that $y \not \in N^+(x)$ and $z \not\in N^+(x)$, then selecting $S \subseteq U-\{y,z\}$ would be a contradiction.
If there were only one vertex $y$ in $U$ such that $y \not \in N^+(x)$, then selecting $S \in U$ with $y\in S$ would be a contradiction.
Therefore, in any cases, we have reached a contradiction.
Hence $U=V(D)$ and so $D$ is isomorphic to $\overleftrightarrow{K}_{t+\lambda}$.
\end{proof}

\begin{Lem}\label{lem:complete,likingvertex}
	Let $D$ be a $(t,\lambda)$-liking digraph for some integers $(t, \lambda) \neq (2,1)$.
	If there exists a vertex $v$ such that $N^+(v)=V(D)-\{v\}$, then $D$ is the complete digraph on $t+\lambda$ vertices. 
\end{Lem}
\begin{proof}
Suppose that there exists a vertex $v$ such that \[N^+_D(v)=V(D) -\{v\}.\]
Then  $D-v$ is a $(t-1,\lambda)$-liking digraph by Proposition~\ref{prop:jan}.
As long as $v$ is an out-neighbor of every other vertex, i.e. \[N^-_D(v)=V(D)-\{v\},\] 
$D-v$ is a $(t,\lambda-1)$-liking digraph.
Then $D-v$ is isomorphic to $\overleftrightarrow{K}_{t+\lambda-1}$ by Lemma~\ref{lem:complete,twoliking}.
Therefore $D$ is isomorphic to $\overleftrightarrow{K}_{t+\lambda}$.
Thus it is sufficient to show $N^-_D(v)=V(D)-\{v\}$.

To this end, we suppose, to the contrary, that
there is a vertex $w$ such that $w \neq v$ and  $w \not \to v$.
Since $|V(D)|\ge t+\lambda \ge t+1$ by Lemma~\ref{lem:outdegree}, there is a subset $T$ of $V(D)-\{v\}$ such that $|T|=t-1$ and $w\in T$.
Then 
$|CN_D^+(T) | \geq \lambda+1$ by Proposition~\ref{prop:out-degree}.
Since $w \not \to v$,
$v \notin CN_D^+(T)$.
Therefore
\[
CN_D^+(T) \subseteq V(D)- \{v\}=N_D^+(v)
\]
and so \[\left|CN_D^+(T\cup\{v\})\right|=\left|N_D^+(v) \cap CN_D^+(T)\right|=\left|CN_D^+(T)\right|\geq \lambda+1.\]
Thus the vertices in $\{v\}\cup T$ have at least $\lambda+1$ common out-neighbors.
However, $|\{v\} \cup T|=t$, which contradicts the hypothesis that $D$ is a $(t,\lambda)$-liking digraph.
Thus $v$ is an out-neighbor of every other vertex in $D$.
\end{proof}

Now, we are ready to prove Theorem~\ref{thm:main2}.

\begin{proof}[Proof of Theorem~\ref{thm:main2}]
	It is obvious that $(a)$ implies each of $(b)$, $(c)$, $(d)$, and $(e)$.
	By Lemma~\ref{lem:complete,twoliking}, $(b)$ implies $(a)$.
	By Lemma~\ref{lem:complete,outdegree}, $(c)$ implies $(a)$.
	Thus $(a)$, $(b)$, and $(c)$ are equivalent.
	If $(t,\lambda) \neq (2,1)$, then Lemma~\ref{lem:complete,likingvertex} guarantees that $(d)$ implies $(a)$.
	If $t\ge 3$, then Proposition~\ref{prop:diregular} guarantees that $(e)$ implies $(a)$.
	Thus the ``furthermore" part is true.
	This completes the proof.
\end{proof}

\section{Concluding Remarks}

\begin{figure}
\begin{center}
 \begin{tikzpicture}[auto,thick]
       \tikzstyle{player}=[minimum size=5pt,inner sep=0pt,outer sep=1pt,draw,circle]

    \tikzstyle{player1}=[minimum size=2pt,inner sep=0pt,outer sep=0pt,fill,color=black, circle]
    \tikzstyle{source}=[minimum size=5pt,inner sep=0pt,outer sep=0pt,ball color=black, circle]
    \tikzstyle{arc}=[minimum size=5pt,inner sep=1pt,outer sep=1pt, font=\footnotesize]
    \path (90:1.5cm)     node [player]  (v1){};
    \path (40:1.5cm)     node [player]  (v7){};
    \path (-10:1.5cm)     node [player]  (v6){};
    \path (-60:1.5cm)     node [player]  (v5){};
 \path (140:1.5cm)     node [player]  (v2){};
  \path (190:1.5cm)     node [player]  (v3){};
    \path (240:1.5cm)     node [player]  (v4){};

   \draw[black,thick,-stealth] (v1) -> (v3);
   \draw[black,thick,-stealth] (v1) -> (v4);
   \draw[black,thick,-stealth] (v1) -> (v6);
   \draw[black,thick,-stealth] (v1) -> (v7);
   \draw[black,thick,-stealth] (v2) -> (v1);
   \draw[black,thick,-stealth] (v2) -> (v4);
   \draw[black,thick,-stealth] (v2) -> (v5);
   \draw[black,thick,-stealth] (v2) -> (v7);
   \draw[black,thick,-stealth] (v3) -> (v1);
   \draw[black,thick,-stealth] (v3) -> (v2);
   \draw[black,thick,-stealth] (v3) -> (v4);
   \draw[black,thick,-stealth] (v3) -> (v6);
   \draw[black,thick,-stealth] (v4) -> (v2);
   \draw[black,thick,-stealth] (v4) -> (v5); 
   \draw[black,thick,-stealth] (v4) -> (v6);
   \draw[black,thick,-stealth] (v4) -> (v7);
   \draw[black,thick,-stealth] (v5) -> (v1);
   \draw[black,thick,-stealth] (v5) -> (v6);
   \draw[black,thick,-stealth] (v5) -> (v7);
   \draw[black,thick,-stealth] (v6) -> (v1);
   \draw[black,thick,-stealth] (v6) -> (v2);
   \draw[black,thick,-stealth] (v6) -> (v3);
   \draw[black,thick,-stealth] (v6) -> (v7);
   \draw[black,thick,-stealth] (v7) -> (v1);
   \draw[black,thick,-stealth] (v7) -> (v3);
   \draw[black,thick,-stealth] (v7) -> (v5);
   \draw[black,thick,-stealth] (v7) -> (v6);

    \end{tikzpicture}
 \end{center}
 \caption{A $(2,2)$-liking digraph which is not diregular}
\label{fig:3-regular}
\end{figure}
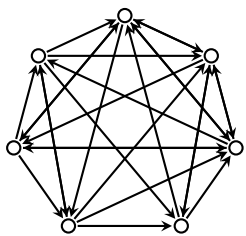

According to \cite{CARSTENS1977286} and \cite{sudolsky1978generalization}, a $(t, \lambda)$-friendship graph is the complete graph on $t+\lambda$ vertices if $t\ge 3$.
This result is partially established for $(t,\lambda)$-liking digraphs with $t\geq \lambda+2$ by our result, Theorem~\ref{thm:main1}.
We would like to know if this digraph version holds in general. 
\begin{Conj}\label{conj2}
	The complete digraph on $t+\lambda$ vertices is the only $(t,\lambda)$-liking digraph if $t\ge 3$.
\end{Conj}
Meanwhile, according to \cite{bose1969graphs},
a $(2, \lambda)$-friendship graph is regular if $\lambda \ge 2$.
Together with the above observation, it is true that a $(t,\lambda)$-friendship graph is regular for any integers $t\geq 2$ and $\lambda \geq 1$ with $(t,\lambda)\neq(2,1)$. 
It is natural to question whether this result also holds for a $(t,\lambda)$-liking digraph for any integers $t\geq 2$ and $\lambda \geq 1$ with $(t,\lambda)\neq(2,1)$.
For the case $t\ge 3$, the answer to the given question immediately follows the equivalence of (a) and (e) in Theorem~\ref{thm:main2}.
\begin{Cor}
If $t\ge 3$, then the complete digraph on $t+\lambda$ vertices is the only diregular $(t,\lambda)$-liking digraph.
\end{Cor}
Yet, this distinction breaks down at $t=2$, as illustrated by our explicit $(2,2)$-liking digraph example that is not diregular (see Figure~\ref{fig:3-regular}).

\section{Acknowledgement}
This work was supported by Science Research Center Program through the National Research Foundation of Korea(NRF) grant funded by the Korean Government (MSIT)(NRF-2022R1A2C1009648).
Hojin Chu was partially supported by a KIAS individual Grant (CG101801) at Korea Institute for advanced Study.

The authors are grateful to Taehee Hong and Homoon Ryu for helping us code the exhaustive search and find a counterexample.


\end{document}